\documentclass[12pt]{elsarticle}

\usepackage{amsmath}
\usepackage{amssymb}
\usepackage{tikz}
\usepackage{verbatim}
\makeatletter
\def\ps@pprintTitle{%
  \let\@oddhead\@empty
  \let\@evenhead\@empty
  \def\@oddfoot{\reset@font\hfil\thepage\hfil}
  \let\@evenfoot\@oddfoot
}
\makeatother

\newtheorem{tetel}{Theorem}
\newtheorem{dfc}[tetel]{Definition}
\newenvironment{biz}{\noindent\textbf{Proof:}\;\relax}{\par\medskip}%
\newtheorem{lemma}[tetel]{Lemma}
\newtheorem{allitas}[tetel]{Claim}
\newtheorem{sejtes}[tetel]{Conjecture}

\newtheorem{cor}[tetel]{Corollary}

\definecolor{Grey}{rgb}{0.5,0.5,0.5}

\begin{document}

\begin{frontmatter}

\title{Bounds on the Number of Edges in Hypertrees}

\author[addr1,addr2]{Gyula Y. Katona\corref{cor1}
}
\ead{kiskat@cs.bme.hu}
\cortext[cor1]{Corresponding author} 
\address[addr1]{Budapest
  University of Technology and Economics\\ Department of Computer
  Science and Information Theory\\
 Magyar Tud\'osok krt. 2, Budapest,  Hungary, H-1117} 
\address[addr2]{HAS Research Group ``Numerical Analysis and Large Networks''\\
P\'azm\'any P\'eter s\'et\'any 1/C, Budapest, Hungary, H-1117 } 

\author[addr1]{P\'eter G. N. Szab\'o
}\ead{szape@math.bme.hu}

\begin{abstract}

Let $\mathcal{H}$ be a $k$-uniform hypergraph. A chain in
$\mathcal{H}$ is a sequence of its vertices such that every $k$
consecutive vertices form an edge. In 1999 Katona and Kierstead
suggested to use chains in hypergraphs as the generalisation of
paths. Although a number of results have been
published on hamiltonian chains in recent years, 
the generalisation of trees with chains has still remained an open area.\\
\indent We generalise the concept of trees for uniform hypergraphs.
We say that a $k$-uniform hypergraph $\mathcal{F}$ is a hypertree 
if every two vertices of $\mathcal{F}$ are connected by a chain, 
and an appropriate kind of cycle-free property holds. An edge-minimal 
hypertree is a hypertree whose edge set is minimal with respect to inclusion.\\
\indent After considering these definitions, we show that a
$k$-uniform hypertree on $n$ vertices has at least $n-(k-1)$ edges up
to a finite number of exceptions, and it has at most $\binom{n}{k-1}$
edges. The latter bound is asymptotically sharp in 3-uniform case.
\end{abstract}

\begin{keyword}
hypertree\sep cycle in hypergraph\sep path in hypergraph

\MSC 05C65 \sep 05C05
\end{keyword}
\end{frontmatter}

\section{Introduction}\label{0}

\subsection{Definitions of Trees}\label{0.1}

A graph is called a \textit{tree} if it is connected and cycle-free (a cycle(-graph) on vertices $v_1, v_2,\ldots, v_n$ is the graph with edges $e_1=\{v_1,v_2\}, e_2=\{v_2,v_3\},\ldots, e_{n-1}=\{v_{n-1},v_n\},$ $e_n=\{v_n,v_1\}$).
Trees play an important role in many areas of applied mathematics,
such as theory of algorithms, cryptography, data
structures and information theory. It is well known that the following are equivalent for a graph $F$ on $n$ vertices.

\begin{itemize}
\item $F$ is a tree,
\item any pair of vertices of $F$ are connected by a unique path,
\item $F$ is an edge-minimal connected graph,
\item $F$ is an edge-maximal cycle-free graph,
\item $F$ is connected and has $n-1$ edges.
\end{itemize}

Let us note that here and in the following, edge-minimal and edge-maximal refers to
minimality or maximality with respect to inclusion of the set of edges.

Our main goal is to generalize the tree propery for hypergraphs. It
turns out that the situation is much more complicated than the graph
case, one can only give lower and upper bounds for the number of edges in hypertrees.

\subsection{Early generalisations of paths and cycles}\label{0.2}

Some earlier concepts of paths and cycles in hypergraphs are recalled
in this subsection. Amongst those, Berge-paths and Berge-cycles are
treated in more details.

A \textit{hypergraph} is a set system i.e. a pair $(V,\mathcal{E})$
where $V$ is a nonempty finite set (called the set of vertices) and
$\mathcal{E}$ is a family of subsets of $V$ (called the set of edges). An
\textit{isolated vertex} in a hypergraph is a vertex which does not
contained in any edges. A \textit{loop} is an edge of cardinality one.

The path-concept, introduced by Berge is one of the very earliest ones
\cite{berge}. We will only consider uniform hypergraphs, but this
definition works for arbitrary hypergraphs, too.

\begin{dfc}\label{Dbergepath}
A Berge-path of length $l$ in a hypergraph $\mathcal{H}$ is a sequence\\
$(v_1,e_1,v_2,e_2,\ldots, v_l,e_l,v_{l+1})$ such that
\begin{itemize}
\item $v_1,v_2,\ldots, v_l, v_{l+1}$ are distinct vertices of $\mathcal{H}$;
\item $e_1,e_2,\ldots,e_l$ are distinct edges of $\mathcal{H}$;
\item $v_k,v_{k+1}\in e_k$, for all $1\leq k\leq l$.
\end{itemize}
If $l>1$ and $v_1=v_{l+1}$, then this ``path'' is called a Berge-cycle of length $l$.
\end{dfc}

\noindent \textbf{Remark:} For graphs the Berge-path and
Berge-cycle is the same as the ordinary path and cycle.
If $\mathcal{H}$ consists of large edges, then we have a considerable freedom in
constructing a path on a given sequence of vertices because every
edge of the path has two fixed vertices (on the path), but the
other vertices can be chosen freely. Our definition
of chain will be more restrictive.

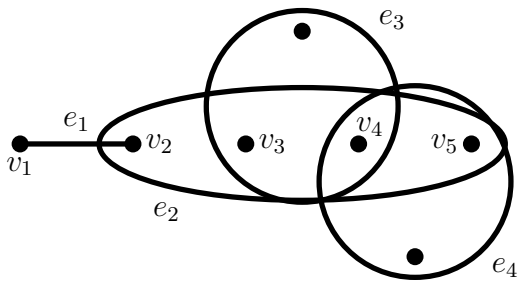
\begin{figure}[h]
\centering
\begin{tikzpicture}[scale=1.5, line width=2pt,
pont/.style={circle, fill=black, inner sep=0.8mm}]
\foreach \x in {0,...,4}
\node at (\x,0) [pont] {};
\draw (0,0) node[below] {$v_{1}$};
\draw (1,0) node[anchor=west] {$v_{2}$};
\draw (2,0) node[anchor=west] {$v_{3}$};
\draw (3.1,0.15) node {$v_{4}$};
\draw (4,0) node[anchor=east] {$v_{5}$};
\node at (2.5,1) [pont] {};
\node at (3.5,-1) [pont] {};
\draw (0,0)--(1,0) node[above,midway] {$e_{1}$};
\draw (2.5,0) ellipse (1.8cm and 0.5cm);
\draw (2.5,1/3) circle (0.85cm);
\draw (3.5,-1/3) circle (0.85cm);
\draw (1.3,-0.6) node {$e_{2}$};
\draw (3.3,1.1) node {$e_{3}$};
\draw (4.3,-1.1) node {$e_{4}$};
\end{tikzpicture}
\caption{A Berge-path on 5 vertices}
\end{figure}

For $u,v\in V$ let $u\equiv v$ denote the fact that there exists a
Berge-path in $\mathcal{H}$ with endpoints $u$ and $v$. It is easy to see
that ``$\equiv$'' is an equivalence relation on $V$ (in contrast to our
chain concept; see Section 2). We define the connected components of
$\mathcal{H}$ to be the equivalence classes of this relation.

We mention two results on Berge-cycle-free hypergraphs. In these,
one can notice analogies with trees and forests of usual graph theory.

\begin{tetel}[Berge \cite{berge}]
  If $\mathcal{H}$ is a hypergraph with $n$ vertices, $p$ connected
  components and $\mathcal{E}=\{e_i\}_{i=1}^m$, then it contains no
  Berge-cycles if and only if $$\sum_{i=1}^m (|e_i|-1)=n-p.$$
\end{tetel}

For graphs this gives that $m=n-p$, which is a well known feature
of forests with $p$ connected components.
 
\begin{tetel}[Lov\'asz \cite{berge}]\label{Tlovász} Let $\mathcal{H}$
be a hypergraph with $n$ vertices $m$ edges and $p$ connected
components, which contains no Berge-cycles of length at least $3$ and no
loops. If $|e_i\cap e_j|\leq 2$ for all $e_i\neq e_j$,
then $$\sum_{i=1}^m (|e_i|-2)<n-p.$$
\end{tetel}

This theorem shows that if $\mathcal{H}$ is a connected $3$-uniform
hypergraph which contains no Berge-cycles of length greater than 2,
then $m<n-1$.

One of the earliest definitions of Hamiltonian cycle in hypergraphs was
given by Bermond et.~al.~in 1976. The authors gave extensions of known
results about Hamiltonian cycles of graphs to the hypergraph setting
\cite{hamchains}.

\begin{dfc} A cyclic permutation $(v_1,v_2, \ldots, v_n)$ of the
vertices of $\mathcal{H}$ is called a hypergraph hamiltonian cycle
if for every $1\leq i\leq n$ there exists an edge $e_i\in\mathcal{E}$
such that $v_i,v_{i+1}\in e_i$.
\end{dfc}

This definition is not equivalent to that of the hamiltonian
Berge-cycle because in this case there can be identical edges among
the $e_i$s.

\subsection{Hamiltonian cycles and Dirac-type theorems for
hypergraphs}\label{0.3}

Although our paper presents a new interesting definition of
hypertrees, it has been inspired by some earlier papers on Dirac-type
theorems for hypergraphs.

The first article in this area is dated to 1999 co-authored by Katona
and Kierstead. It remains an active field of research, and a number of
results have been published on hamiltonian chains in recent years by
E. Szemer\'edi, V. R\"odl, A. Ruci\'nski, D. K\"uhn, D. Osthus, R. Mycroft,
H. H\`an, M. Schacht and others.
For more details  see
\cite{looseham,lcycles,hamchains,loosecyc,approx, szemeredi}.


\section{Definition of hypertrees}\label{1}

In this section, we generalise the concept of trees for $k$-uniform
hypergraphs, where chains   play the role of paths. After
introducing the basic definitions, we discuss lower and upper
bounds for the edge number of hypertrees and show that these are sharp in particular cases.

\subsection{Possibilities for defining hypertrees}\label{1.1}

First, one has to clarify the notions of cycle, semicycle and chain.
The relation between a chain and a path is similar to the relation between
a tight hamiltonian cycle (see \cite{loosecyc}) and a usual hamiltonian cycle.  In the
following sections we assume that there are no multiple edges in a
hypergraph.

\begin{dfc}[Cycle] The $k$-uniform hypergraph
$\mathcal{C}=(V,\mathcal{E})$ is a cycle if there exists a cyclic
sequence $v_1,v_2,\ldots, v_l$ of its vertices such that every vertex
appears at least once (possibly more times) and for all $1\leq i\leq l$,
$\{v_{i}, v_{i+1},\ldots, v_{i+k-1}\}$ are distinct edges of $\mathcal{C}$. The length of the cycle
$\mathcal{C}$ is the number of its edges: $l$.
\end{dfc}

Note that a cycle has at least $k+1$ edges. To see this, notice that if
it is defined by the sequence $v_1,v_2,\ldots, v_l$, then it has
exactly $l$ (different) edges by definition.  $k$-uniformity implies
that $l$ is at least $k$. However, $l=k$ means that the $k$ edges of
the definition coincide, each covers the whole vertex set.

Cycles are too special structures for our purpose, so we use a weaker
concept. If a chain intersects itself, then it contains a
subhypergraph called a semicycle. Hypergraphs without semicycles are
in close resemblance with (ordinary) forests.

\begin{dfc}[Semicycle] The nonempty $k$-uniform hypergraph
$\mathcal{C}=(V,\mathcal{E})$ is a semicycle if there exists a
sequence $v_1,v_2,\ldots, v_l$ of its vertices such that every vertex
appears at least once (possibly more times), $v_1=v_l$ and
for all $1\leq i\leq l-k+1$, $\{v_{i}, v_{i+1},\ldots, v_{i+k-1}\}$ 
are distinct edges of $\mathcal{C}$. The length of the semicycle
$\mathcal{C}$ is the number of its edges: $l-k+1$.
\end{dfc}

Notice that a semicycle must have at least $3$ edges. Clearly, one edge
cannot form a semicycle. If a semicycle has only two edges, then it is
defined by a sequence $v_1, v_2,\ldots, v_k,v_{k+1}$ of $k+1$
vertices. However, $v_1=v_{k+1}$, so there are $k$ different vertices
for two different edges, which is not enough.

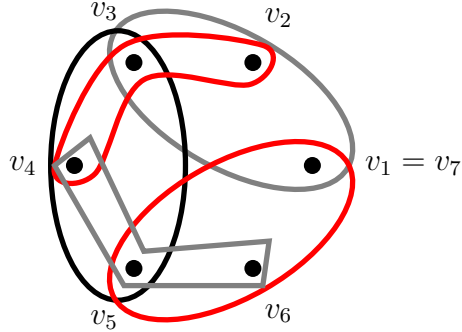
\begin{figure}[h] \centering
\begin{tikzpicture} [scale=1.5, line width=2pt, pont/.style={circle, fill=black, inner
sep=0.8mm}]

\node at (0:30pt) [pont] {}; \node at (60:30pt) [pont] {}; \node at
(120:30pt) [pont] {}; \node at (180:30pt) [pont] {}; \node at
(240:30pt) [pont] {}; \node at (300:30pt) [pont] {};

\draw[Grey, rotate=60*1] (0:19pt) ellipse (17pt and 34pt);
\draw[rotate=60*3] (0:19pt) ellipse (17pt and 34pt);
\draw[red, rotate=60*5] (0:19pt) ellipse (17pt and 34pt);

\draw[red] plot[smooth cycle, tension=0.6]
coordinates{(60*2-65:25pt) (60*2-60:35pt) (60*2:35pt)
(60*2+60:35pt) (60*2+65:25pt) (60*2:25pt)};

\draw[Grey] plot coordinates{(60*4-75:27pt) (60*4-60:35pt) (60*4:35pt)
(60*4+60:35pt) (60*4+75:27pt) (60*4:25pt)}--cycle;

\foreach \x in {2,4,6} \draw (60*\x-60:43pt) node {$v_{\x}$};
\foreach \x in {3,5} \draw (60*\x-60:45pt) node {$v_{\x}$};
\draw (0:40pt) node[right] {$v_1=v_7$};
\end{tikzpicture}
\caption{A non-self-intersecting $3$-uniform semicycle of length 5}
\end{figure}

Simple chains will play the most important role in defining
hypertrees because we intend to require a natural chain-connectedness
property.

\begin{dfc}[Chain] The nonempty $k$-uniform hypergraph
$\mathcal{L}=(V,\mathcal{E})$ is a chain if there exists a sequence
$v_1,v_2,\ldots, v_l$ of its vertices such that every vertex appears
at least once (possibly more times), $v_1\neq v_l$ and 
for all $1\leq i\leq l-k+1$, $\{v_{i}, v_{i+1},\ldots, v_{i+k-1}\}$ 
are distinct edges of $\mathcal{C}$.  The length of the chain
$\mathcal{L}$ is the number of its edges: $l-k+1$.
\end{dfc}

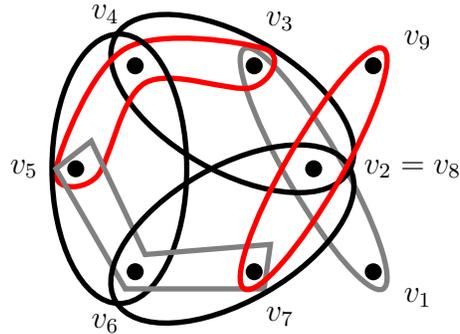
\begin{figure}[h] \centering
\begin{tikzpicture} [scale=1.5, line width=2pt, pont/.style={circle, fill=black, inner
sep=0.8mm}]

\node at (0:30pt) [pont] {}; \node at (60:30pt) [pont] {}; \node at
(120:30pt) [pont] {}; \node at (180:30pt) [pont] {}; \node at
(240:30pt) [pont] {}; \node at (300:30pt) [pont] {}; \node at
(330:52pt) [pont] {}; \node at (30:52pt) [pont] {};

\draw[Grey, rotate=300] (60:30pt) ellipse (35pt and 7pt);
\draw[rotate=60*1] (0:19pt) ellipse (17pt and 34pt);
\draw[red] plot[smooth cycle, tension=0.6]
coordinates{(60*2-65:25pt) (60*2-60:35pt) (60*2:35pt)
(60*2+60:35pt) (60*2+65:25pt) (60*2:25pt)};
\draw[rotate=60*3] (0:19pt) ellipse (17pt and 34pt);
\draw[Grey] plot coordinates{(60*4-75:27pt) (60*4-60:35pt) (60*4:35pt)
(60*4+60:35pt) (60*4+75:27pt) (60*4:25pt)}--cycle;
\draw[rotate=60*5] (0:19pt) ellipse (17pt and 34pt);
\draw[red, rotate=60] (300:30pt) ellipse (35pt and 7pt);

\draw (330:65pt) node {$v_1$};
\draw (0:55pt) node {$v_{2}=v_8$};
\foreach \x in {4,6} \draw (60*\x-120:45pt) node {$v_{\x}$};
\foreach \x in {3,5,7} \draw (60*\x-120:43pt) node {$v_{\x}$};
\draw (30:65pt) node {$v_9$};

\end{tikzpicture}
\caption{A self-intersecting $3$-uniform chain of length 7}
\end{figure}

A chain (cycle) is \textit{self-intersecting} if there is a vertex
appearing at least twice in its defining sequence.  Similarly, a
semicycle is \textit{self-intersecting} if there is a vertex appearing
at least twice in its defining sequence, except for the condition that
the first and the last vertices must be identical.

\begin{figure}[h] \centering
\begin{tikzpicture} [scale=1.5, line width=2pt,pont/.style={circle, fill=black, inner
sep=0.8mm}]

\node at (0,0) [pont] {}; \node at (1,0) [pont] {}; \node at (2,0)
[pont] {}; \node at (3,0) [pont] {}; \node at (4,0) [pont] {}; \draw
(1,0) ellipse (35pt and 8pt); \draw (2,0) ellipse (35pt and 8pt);
\draw (3,0) ellipse (35pt and 8pt);

\end{tikzpicture}
\caption{A non-self-intersecting $3$-uniform chain of length 3}
\end{figure}
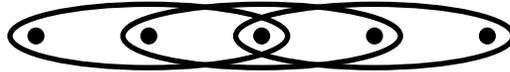
\noindent \textbf{Remark:} A non-self-intersecting chain (cycle) is also called a
tight-path (tight-cycle) in other related papers.

\noindent \textbf{Remark:} The length of a non-self-intersecting chain
on $n$ vertices is $n-(k-1)$, matching the fact that the length of a ($2$-uniform)
path on $n$ vertices is $n-1$.

 Similarly, the length of a
non-self-intersecting semicycle on $n$ vertices is $n-k+2$, and for
graphs, the length of a (2-uniform) cycle on $n$ vertices is $n$.

\begin{dfc} A $k$-uniform hypergraph $\mathcal{H}$ is
\begin{itemize}
\item \emph{chain-connected} if every pair of its vertices is
connected by a chain, i.e. there exists a subhypergraph which
is a chain and contains both vertices;
\item \emph{semicycle-free} if it contains no semicycle as 
subhypergraph;
\item \emph{edge minimal/maximal} with respect to the property $\phi$
if it has property $\phi$, and its edge set is minimal/maximal with
respect to inclusion.
\end{itemize}
\end{dfc}

The following claim shows that to verify that a hypergraph is
semicycle-free it is enough to show that it does not contain any
non-self-intersecting semicycles.

\begin{allitas}\label{Tsemicyc} If a $k$-uniform hypergraph
$\mathcal{H}$ contains a semicycle, then it also contains a
non-self-intersecting semicycle.
\end{allitas}

\begin{biz}
  Let $C$ be the shortest semicycle in $\mathcal{H}$, defined by the
  vertex sequence $(x_1, x_2, \ldots, x_l)$. We show that $C$ is
  non-self-intersecting. Assume on the contrary that $C$ intersects
  itself. Then there exist indices $1\leq i<j\leq l$ such that
  $x_i=x_j$ and $i\neq 1$ or $j\neq l$. However, the sequence $x_i,
  x_{i+1},\ldots, x_j$ defines a semicycle which is shorter than $C$,
  a contradiction. $\square$
\end{biz}

\begin{cor} If a uniform hypergraph does not contain any non-self-intersecting
semicycles, then it is semicycle-free.
\end{cor}

The counterpart of the previous theorem among graphs states that if a
graph contains a closed walk, then it also contains a cycle. Among
hypergraphs, this theorem becomes false if we use 
cycle instead of semicycle. So, considering the semicycle as the
generalisation of the graph-theoretical cycle seems to be a fruitful
idea.

\begin{allitas}
A semicycle-free $k$-uniform hypergraph does not contain a self-intersecting chain.
\end{allitas}

\begin{biz} Let $L=(x_1,\ldots,x_l)$ be a chain in $\mathcal{H}$,
and assume indirectly that it is self-intersecting. Then there exists
$1\leq i<j\leq l$ such that $x_i=x_j$ ($i\neq1$ or $j\neq l$ since
$x_1\neq x_l$). Now the part $C=(x_i, \ldots, x_j)$ of the chain induces a
semicycle, which is a contradiction. $\square$
\end{biz}

At this point, we can define hypertrees the same way we defined trees
in Section \ref{0}. Some of these definitions are not compatible with
the concept of chain, while others may be too general. One has to take
into consideration that chain-connectedness is not a transitive
property, and that there are two different possibilities to generalize cycles.

 Let
$\mathcal{F}$ be a $k$-uniform hypergraph. Let us explore the various possible 
generalisations of the equivalent tree definitions.

\begin{enumerate}
\item $F$ is connected and cycle-free: here there are two
possibilities: $\mathcal{F}$ is chain-connected and cycle-free, or
$\mathcal{F}$ is chain-connected and semicycle-free.

We require the stronger semicycle-freeness condition instead of
cycle-freeness for hypertrees. Although cycle-freenes is a weaker
condition, later it turns out that it is not that much weaker, but
there will be a proof where the
absence of semicycles will play an important role.

\item Every pair of vertices of $F$ is connected by exactly one path:
  this is not a useful definition in this context, because in this case
  there could be no chains in $\mathcal{F}$ with length more than
  one: If two edges $e,f\in\mathcal{E}$ intersect each other
  in $2$ vertices, for example, $e\cap f=\{u,v\}$, then $u$ and $v$
  are connected by two chains ($e$ and $f$) of length one, 
  contradicting the definition.

 So, $\mathcal{F}$ can only be a kind of block-design, which are
already well investigated,
therefore we rejected
this idea.

\item $F$ is an edge-minimal connected graph: the generalisation of
  this definition is that $\mathcal{F}$ is an edge-minimal
  chain-connected hypergraph. This is not a subcase of the first
  definition unless we require the semicycle-free property. The
  example below points to following counterexample proves this
  observation.

Let $n\geq 5$, $V=\{x_1, x_2, \ldots, x_n, y_1, y_2, \ldots, y_n, z_1,
z_2, \ldots, z_n\}$ and $\mathcal{C}$ be the $3$-uniform cycle
determined by the vertex-sequence $x_1,x_2,\ldots,x_n$. Furthermore,
let \[\mathcal{E}_y=\{\{x_{i-1},x_i,y_i\}\colon 1\leq i\leq n\},\quad
\mathcal{E}_z=\{\{z_i,x_i,x_{i+1}\}\colon 1\leq i\leq n\}\] (the indices are
understood cyclicly). Then
$\mathcal{F}=(V,\mathcal{E}(\mathcal{C})\cup\mathcal{E}_y\cup\mathcal{E}_z)$
is edge-minimal chain-connected hypergraph, but it is not
semicycle-free, moreover it is not cycle-free because $\mathcal{C}$
is its subhypergraph.

\item $F$ is edge-maximal cycle-free graph: the generalisation of this
definition is that $\mathcal{F}$ is an edge-maximal semicycle-free
hypergraph. Again, this is not a subcase of the first definition
unless we require chain-connectedness. Instead of semicycle-freeness
one could also require cycle-freeness here.
\end{enumerate}

To sum these ideas, we would like to give definitions which are strong
enough, but yield us not only the trivial structures.

\begin{dfc}[Hypertree]\label{Dhypert} The $k$-uniform hypergraph
$\mathcal{F}$ is a hypertree if it is chain-connected and
semicycle-free.
\end{dfc}

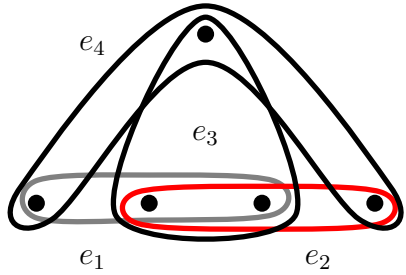
\begin{figure}[h] \centering
\begin{tikzpicture} [scale=1.5, line width=2pt,pont/.style={circle, fill=black, inner
sep=0.8mm}]

\node at (0,0) [pont] {}; \node at (1,0) [pont] {}; \node at (2,0)
[pont] {}; \node at (3,0) [pont] {}; \node at (1.5,1.5) [pont] {};

\draw[Grey] plot[smooth cycle, tension=0.6] coordinates{(0,-0.15+0.02)
(0,0.15+0.05) (1,0.2+0.05) (2.1,0.15+0.05) (2.1,-0.15+0.05)
(1,-0.2+0.05)};
\draw[red] plot[smooth cycle, tension=0.6]
coordinates{(0.9,-0.15-0.03) (0.9,0.15-0.05) (2,0.2-0.05)
(3.05,0.15-0.05) (3.05,-0.15-0.03) (2,-0.2-0.03)};
\draw plot[smooth cycle, tension=0.6] coordinates{(0.141,-0.141) (-0.141,0.141)
(1.5,1.75) (3.141,0.141) (2.859,-0.141) (1.5,1.25)};
\draw plot[smooth cycle, tension=0.6] coordinates{(1-0.3,0-0.1) (1.5,1.65)
(2+0.3,0-0.1)};

\draw (0.5,-0.5) node { $e_1$}; \draw (2.5,-0.5) node
{ $e_2$}; \draw (1.5,0.6) node { $e_3$};
\draw (0.5,1.4) node { $e_4$};

\end{tikzpicture}
\caption{A 3-uniform hypertree of order 5}\label{Ahipfa}
\end{figure}

\begin{dfc}[Edge-minimal hypertree]\label{Dedgemin} $\mathcal{F}$ is
an edge-minimal hypertree if it is a hypertree, and deleting any of its edge
$e$, $\mathcal{F}\backslash\{e\}$ is not a hypertree any more
(i.e. chain-connectedness does not hold).
\end{dfc}

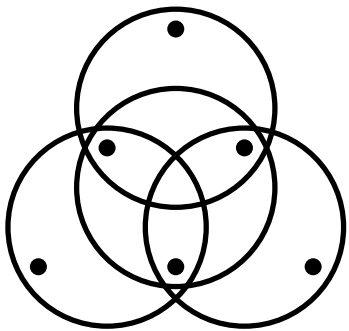
\begin{figure}[h] \centering
\begin{tikzpicture} [scale=1.5, line width=2pt,pont/.style={circle, fill=black, inner
sep=0.8mm}]

\node at (-90:20pt) [pont] {}; \node at (30:20pt) [pont] {}; \node at
(150:20pt) [pont] {}; \node at (-30:40pt) [pont] {}; \node at
(90:40pt) [pont] {}; \node at (210:40pt) [pont] {};

\draw (0:0pt) circle (25pt); \draw (-30:20pt) circle (25pt); \draw
(90:20pt) circle (25pt); \draw (210:20pt) circle (25pt);

\end{tikzpicture}
\caption{An edge-minimal $3$-uniform hypertree on 6 vertices}
\end{figure}

\begin{dfc}[Edge-maximal hypertree]\label{Dedgemax} $\mathcal{F}$ is
an edge-maximal hypertree if it is a hypertree, and adding any new
edge $e$ to it, $\mathcal{F}\cup\{e\}$ is not a hypertree any more
(i.e. semicycle-freeness does not hold).
\end{dfc}

\begin{figure}[h] \centering
\begin{tikzpicture} [scale=1.5, line width=2pt, pont/.style={circle, fill=black, inner
sep=0.8mm}]

\node at (0,0) [pont] {}; \node at (1,0) [pont] {}; \node at (2,0)
[pont] {}; \node at (-1,0) [pont] {}; \node at (0.5,1.5) [pont] {};
\node at (0.5,-1.5) [pont] {};

\draw[Grey] plot[smooth cycle, tension=0.6] coordinates{(-1,-0.15+0.02)
(-1,0.15+0.05) (0,0.2+0.05) (1.08,0.15+0.09) (1.08,-0.15+0.05)
(0,-0.2+0.05)};
\draw[red] plot[smooth cycle, tension=0.6]
coordinates{(-0.1,-0.15-0.03) (-0.1,0.15-0.05) (1,0.2-0.05)
(2.05,0.15-0.05) (2.05,-0.15-0.03) (1,-0.2-0.03)};
\draw plot[smooth cycle, tension=0.6] coordinates{(-0.86,-0.14) (-1.14,0.14)
(0.5,1.75) (2.141,0.141) (1.859,-0.141) (0.5,1.25)};
\draw[Grey] plot[smooth cycle, tension=0.6] coordinates{(-0.3,0-0.1) (0.5,1.65)
(1+0.3,0-0.1) (0.5,-0.5)};
\draw[red] plot[smooth
cycle, tension=0.55] coordinates{(-1+0.106,0+0.3) (-1-0.3,0-0.2)
(0.5,-1.5-0.4) (2+0.3,0-0.2) (2-0.106,0+0.3) (0.5,-1.5+0.3)};
\draw plot[smooth cycle, tension=0.7]
coordinates{(-0.14,0.0) (0.5,0.6) (1.14,0.0) (0.5,-1.7)};

\end{tikzpicture}
\caption{An edge-maximal $3$-uniform hypertree on 6 vertices}
\end{figure}
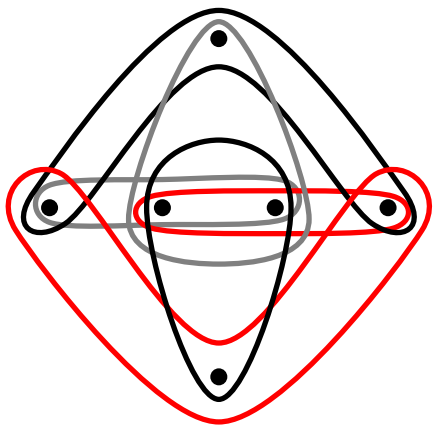

In this way, the edge maximal/minimal hypertrees are also 
hypertrees, but they are  the extreme cases among hypertrees.

The main reason to use semicycle-free
property is that Claim \ref{Tsemicyc} implies that every chain is
non-self-intersecting in a hypertree. Without this property we should
face with substantially more complicated case-analysis.

We mentioned above that being ``connected by chain'' is not a
transitive property, thus it is not an equivalence relation. This
fact is the main difficulty in this topic. It also means that it is usually
hard to mimic the graph theoretic proofs, we need to invent new ideas.

\begin{dfc}[$l$-hypertree] A hypertree is called an $l$-hypertree if
all chains in it have length at most $l$.
\end{dfc}

In Section \ref{4} it will be shown that there are better bounds for the
edge number of an $l$-hypertree if $l<k$.

To every hypergraph we assign a graph, which more or less preserves
the structure of chains. This will be a useful tool in some proofs.

\begin{dfc}[Tight line graph]\label{Dedgegraph} Let
$\mathcal{H}=(V,\mathcal{E})$ be a $k$-uniform hypergraph. The tight
line graph of $\mathcal{H}$ is the graph
$L_{\mathcal{H}}=(\mathcal{E},E)$ where
$E=\{\{e,f\}\colon e,f\in\mathcal{E}, |e\cap f|=k-1\}$.
\end{dfc}
This definition is not equivalent to the definition of the
usual line graph, where all pairs of intersecting edges are
adjacent. In our case we are only interested in the substantial
intersections (of size $(k-1)$) since the edges of chains are joined
in this way.
\begin{dfc} The hypergraph $\mathcal{H}$ is line-graph-connected if
its tight line graph is connected.
\end{dfc}
 This will be used in the proof of the lower bound for
the edge number of hypertrees. The proof will be easier for
line-graph-connected hypertrees, while the general case can be traced
back to the tight line graph components.

\begin{figure}[h] \centering
\begin{tikzpicture} [scale=1.5, line width=2pt,pont/.style={circle, fill=black, inner
sep=0.8mm}]

\node at (0,0) [pont] {}; \node at (0,1) [pont] {}; \node at
(0.866,-0.5) [pont] {}; \node at (-0.866,-0.5) [pont] {};

\draw (0,0) node[anchor=south] { $e_3$}; \draw (0,1)
node[anchor=south] { $e_4$}; \draw (0.866,-0.5)
node[anchor=north west] { $e_2$}; \draw (-0.866,-0.5)
node[anchor=north east] { $e_1$};

\draw (0,0)--(0.866,-0.5); \draw (0,0)--(-0.866,-0.5); \draw
(0.866,-0.5)--(-0.866,-0.5);

\end{tikzpicture}
\caption{The tight line graph of the hypertree of
figure \ref{Ahipfa}, which is not
line-graph-connected.}
\end{figure}
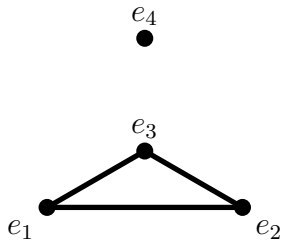

\subsection{Basic types of edge-minimal hypertrees}\label{1.2}

In the following subsection we are going to review some special types
of edge-minimal hypertrees. We will sometimes refer back to these
constructions later.

\begin{dfc}[Tight star]\label{Dstar} The $k$-uniform hypergraph
$\mathcal{S}_n$ of order $n$ is a tight star if $n\geq k$ and there exist
$u_1, u_2, \ldots, u_{k-1}\in V(\mathcal{S}_n)$ such that
\[\mathcal{E}(\mathcal{S}_n)=\{\{u_1,u_2, \ldots, u_{k-1} ,w\}\colon  w\in
V(\mathcal{S}_n), w\neq u_i, \text{ for } 1\leq i\leq k-1\}.\]
\end{dfc}
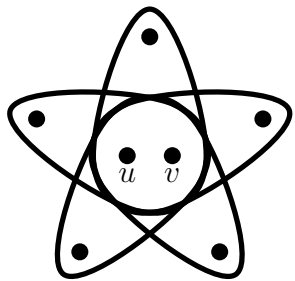
\begin{figure}[h] \centering
\begin{tikzpicture} [scale=1.5, line width=2pt,pont/.style={circle, fill=black, inner
sep=0.8mm}]

\node at (-0.2,0) [pont] {}; \node at (0.2,0) [pont] {}; \draw
(-0.2,0) node[anchor=north] {$u$}; \draw (0.2,0)
node[anchor=north] {$v$}; \foreach \x in {0,...,4} \node
at (90+72*\x:30pt) [pont] {};

\foreach \x in {0,...,4} 
\draw[rotate=72*\x] plot[smooth cycle, tension=0.8]
coordinates{(-0.5,0) (0,1.3) (0.5,0) (0,-0.5)};
\end{tikzpicture}
\caption{The $3$-uniform tight star with 7 vertices}
\end{figure}

\begin{allitas}
Every tight star is an edge-minimal hypertree.
\end{allitas}

Note, that the two most simple class of hypertrees are the tight stars
and the non-self-intersecting chains (tight paths). From now on, we
will write star instead of tight star.

Non-self-intersecting chains and semicycles will be denoted by roman
capital letters, usually by $L$ and $C$.

\begin{dfc}[$l$-Flower]\label{Dlflower} The hypergraph $\mathcal{V}_n$
of order $n$ is an $l$-flower if $n> k>l\ge 1$,
$V(\mathcal{V}_n)=\{v_1,\ldots, v_{n-l}, u_1, \ldots, u_l\}$ and
\[\mathcal{E}(\mathcal{V}_n)=\{\{v_i, v_{i+1}, \ldots, v_{i+k-l-1},u_1,
\ldots, u_l\}\colon 1\leq i\leq n-l\},\] where  indices are understood cyclicly.  
A $1$-flower is simply called flower.
\end{dfc}

Note, that a $(k-1)$-flower is the same as a star.

\begin{figure}[h] \centering
\begin{tikzpicture} [scale=1.5, line width=2pt,pont/.style={circle, fill=black, inner
sep=0.8mm}]

\node at (0,0) [pont] {}; \foreach \x in {0,...,4} \node at
(90+72*\x:30pt) [pont] {};

\foreach \x in {0,...,4} 
\draw[rotate=72*\x] (126:18pt) circle (22pt);
\end{tikzpicture}
\caption{The 3-uniform flower with 6 vertices}
\end{figure}
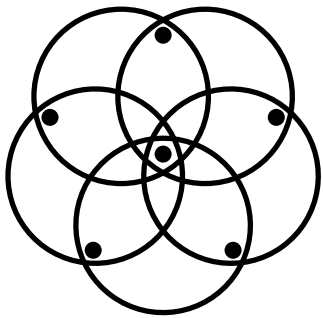

\begin{allitas} Let $2<k<n$. The $k$-uniform flower $\mathcal{V}_n$ is
  a hypertree if and only if $2k-1\leq n\leq 4(k-1)$. It is an 
   edge-minimal hypertree if and only if $3(k-1)\leq n\leq 4(k-1)$.
\end{allitas}

The proof is quite simple and left to the reader. The key is the fact,
that $v_i$ and $v_j$ are connected by a chain in $\mathcal{V}_n$ if
and only if $\min\{|j-i|,n-1-|j-i|\}<2k-2$. $v_i$ and $v_j$ are
connected by two edge-disjoint chains (or $v_i=v_j$ if it is the first
vertex of a semicycle) if and only if $\max\{|j-i|,n-1-|j-i|\}<2k-2$.

\begin{dfc}[Focus-vertex] Let $\mathcal{F}=(V,\mathcal{E})$ be an
edge-minimal hypertree. Then $v\in V$ is a focus-vertex of
$\mathcal{F}$ if it is contained in every edge of $\mathcal{F}$.
\end{dfc}

Obviously, an edge-minimal hypertree may have at most $k$
focus-vertices, and only the trivial hypertree (which consists of 1
edge) has exactly $k$ of them. Hypertrees with $k-1$ focus-vertices
are stars, and flowers have one focus-vertex.
It is an interesting question in itself that what is the maximum
number of edges in a 3-uniform hypertree on $n$ vertices with exactly
1 focus-vertex. Removing this focus-vertex leads to a
graph-theoretical problem.

\begin{dfc}[Geometric hypertree] A hypertree $\mathcal{F}$ is
$l$-geometric if every $l$-element set of its vertices is contained in
exactly one edge.
\end{dfc}

For example, the coincidence-hypergraph of the Fano-plane and other
Steiner-systems are geometric hypertrees. $l$-geometric hypertrees are
exactly the $l$-$(n,k,1)$ block-designs. Therefore, every geometric
hypertree with the same parameters has an equal number of edges that
can be computed from the parameters. (We will see that this is not
true for all hypertrees.)

\begin{itemize}
\item Every $k$-uniform $2$-geometric hypertree is edge-minimal.
\item A $k$-uniform $l$-geometric hypertree has exactly
  $\frac{1}{\binom{k}{l}} \binom{n}{l}$ edges which is well known from
  the theory of block-designs.
\end{itemize}

\begin{dfc}[Recursive hypertree]\label{Drek} An edge-minimal hypertree
$\mathcal{F}$ is recursive if it can be obtained by the following recursive
construction:
\begin{enumerate}
\item One edge consisting of $k$ vertices is a recursive
hypertree.
\item If we add a new vertex $v$ and some new edges all containing $v$ to a
recursive hypertree such that the resulting hypergraph is still an edge-minimal hypertree,
then this is also a recursive hypertree.
\end{enumerate}
\end{dfc}

Can we build up all edge-minimal hypertrees with recursive construction? The answer
is negative: deleting any vertex $v_i$ and the two edges containing
it from the 3-uniform flower $\mathcal{V}_7$, we get a
non-chain-connected hypergraph (see Definition \ref{Dlflower}) since
$v_{i-1}$ and $v_{i+1}$ become chain-disconnected. If we delete the
focus-vertex $u_1$ in turn, then we have to erase all edges, so the
remaining hypergraph obviously not chain-connected. This shows that
$\mathcal{V}_7$ is not a recursive hypertree.

\section{Lower bound for the edge number of hypertrees}\label{2}

A tree on $n$ vertices has exactly $n-1$ edges. In case of hypertrees
the situation is more complicated. Chains and stars have $n-(k-1)$
edges, however,  $(k-1)$-geometric hypertrees have $\frac{1}{k}
\binom{n}{k-1}$ edges. Anyway, it seems like $n-(k-1)$ is the tight
lower bound for the number of edges, however, it turns out that it is
true only if $n>n_0(k)$.

\begin{tetel}\label{Tlower} Let $\mathcal{F}$ be a $k$-uniform
chain-connected hypergraph with $n$ vertices and $m$ edges. If $n\geq
(k-1)^2$, then $m\geq n-(k-1)$.
\end{tetel}

\begin{biz} Let $L$ be the tight line graph of $\mathcal{F}$, and
let $\mathcal{E}_1, \mathcal{E}_2, \ldots, \mathcal{E}_l$ denote the
vertex sets of the connected components of $L$. Furthermore, let
$V_i=\bigcup \{e\colon e\in\mathcal{E}_i\}$, $i=1,\ldots, l$ be the projections of the
components to the vertex set $V$, called the classes of
$\mathcal{F}$. Then it can be easily shown that the following hold:
\begin{itemize}
\item $\bigcup_{i=1}^l \mathcal{E}_i=\mathcal{E}$,
$\mathcal{E}_i\cap\mathcal{E}_j=\emptyset$ if $i\neq j$.
\item $\bigcup_{i=1}^l V_i=V$ since every vertex of $\mathcal{F}$ is
contained in some edges.
\item $\forall u,v\in V$, $\exists i$ such that $u,v\in V_i$.
\end{itemize} The last point follows from the fact that $\mathcal{F}$
is chain-connected and the edges of a chain form the vertex set of a
path in the tight line graph. Therefore, the edges of the chain that
connects $u$ and $v$, are in the same tight line graph component
$\mathcal{E}_i$ of $L$, whose projection $V_i$ contains the vertices
of the chain, in particular both $u$ and $v$.

 Now we are going to prove
a lemma which says that the theorem is true for line-graph-connected
hypertrees, so it is true for every class.

\begin{lemma} If $\mathcal{H}=(V',\mathcal{E}')$ is a
line-graph-connected $k$-uniform hypergraph that contains no isolated
vertex, then $|\mathcal{E}'|\geq |V'|-(k-1)$.
\end{lemma}

\begin{biz} Let $e\in\mathcal{E}'$ be an arbitrary edge and
$\mathcal{X}_1=(e,\{e\})$ be a subhypergraph. We will define a
subhypergraph-sequence $\mathcal{X}_1 \subset \mathcal{X}_2
\subset \ldots \subset \mathcal{X}_m$, for which
$\mathcal{X}_m=\mathcal{H}$ and for every $i$,
$i=|\mathcal{E}(\mathcal{X}_i)|\geq |V(\mathcal{X}_i)|-(k-1)$. This
implies the lemma.

 Let $m=|\mathcal{E}'|$. Assume that we
have already defined $\mathcal{X}_{i-1}$ for some $i-1<m$, and let $L$
denote the tight line graph of $\mathcal{H}$. Since we are not in the
$m\text{th}$ step and $L$ is connected, there must be an edge in $L$
between $\mathcal{E}(\mathcal{X}_{i-1})$ and
$\mathcal{E}'\backslash\mathcal{E}(\mathcal{X}_{i-1})$. Let
$f_{i-1}\in \mathcal{E}(\mathcal{X}_{i-1})$ and $e_i\notin
\mathcal{E}(\mathcal{X}_{i-1})$ denote the two endpoints of this
edge. Then $|f_{i-1}\cap e_i|=k-1$, so $|e_i\cap
V(\mathcal{X}_{i-1})|\geq k-1$. Let
$\mathcal{X}_i=\mathcal{X}_{i-1}\cup \{e_i\}=(V(\mathcal{X}_{i-1})\cup
e_i,\mathcal{E}(\mathcal{X}_{i-1})\cup\{e_i\})$.

 Thus
starting from $\mathcal{X}_1$, in every step we add a new edge, and we continue
this process until all the edges of $\mathcal{H}$ are added. At the last
step $\mathcal{X}_m=\mathcal{H}$ will hold since there is no isolated vertex.

 By induction, it can be shown that for every index $i$,
$|\mathcal{E}(\mathcal{X}_i)|\geq |V(\mathcal{X}_i)|-(k-1)$. If $i=1$,
then the claim is trivially true. Assume that the claim is true for
$i-1$, namely, $|\mathcal{E}(\mathcal{X}_{i-1})|\geq
|V(\mathcal{X}_{i-1})|-(k-1)$. In the $i^{\text{th}}$ step we add
one new edge and at most one new vertex to $\mathcal{X}_{i-1}$
because $e_i$ intersects $V(\mathcal{X}_{i-1})$ in at least $k-1$
vertices. Thus $|\mathcal{E}(\mathcal{X}_i)|\geq
|V(\mathcal{X}_i)|-(k-1)$ also holds. $\square$
\end{biz}

Continuing the proof of the theorem, all of the hypergraphs
$(V_i,\mathcal{E}_i)$ meet the conditions of the lemma, so
$|\mathcal{E}_i|\geq |V_i|-(k-1)$ for $i=1,2,\ldots,
l$. Therefore, $$|\mathcal{E}|=\sum_{i=1}^l |\mathcal{E}_i|\geq
\sum_{i=1}^l \left(|V_i|-(k-1)\right)=\left(\sum_{i=1}^l
|V_i|\right)-l(k-1).$$ Let $\sigma=\sum_{i=1}^l |V_i|$. For proving $m\geq
n-(k-1)$ it is enough to show that $\sigma-l(k-1)\geq n-(k-1)$, or
equivalently $$\sigma\geq n+(l-1)(k-1).$$

Notice that $\sigma\geq n$ certainly holds since every vertex of
$\mathcal{F}$ is covered by one of the classes.

 We may
assume that $|V_i|<n$ for every $i$, otherwise there would be a
line-graph component $(V_j,\mathcal{E}_j)$ such that $|V_j|=n$ and
$n-(k-1)=|V_j|-(k-1)\leq |\mathcal{E}_j|\leq |\mathcal{E}|$, which
would finish the proof.

 We introduce a
parameter $r=\min_{x\in V} |\{i\colon x\in V_i\}|$, which helps us in
splitting the problem into two parts. It is clear that $r\geq 2$
always holds since if $r=1$, then there would be a vertex $x$, which
is covered by only one class, say $V_j$. However, for every $y\in V$
the pair $\{x,y\}$ of vertices has to be covered by a class, which can
only be $V_j$. Hence $V_j=V$, which
contradicts the previous assumption.
\medskip

 To complete the proof,
we only need two inequalities.

\begin{lemma} \mbox{}
\begin{enumerate} \renewcommand{\labelenumi}{\textit{(\theenumi)}}
\item $\sigma\geq rn$,
\item $\sigma\geq n+r-1+(l-r)k$.
\end{enumerate}
\end{lemma}

\begin{biz}

\noindent $(1)$ We are going to compute $\sigma$ in two ways. We can sum up the
size of the classes or we can sum up for all vertices the number of
classes containing them. Thus $\sigma=\sum_{x\in V} |\{i\colon  x\in V_i\}|\geq rn$,
which verifies the first inequality.

\noindent $(2)$ For every vertex $x$ the union of the classes
containing $x$ must be $V$, since they have to cover every pair $\{x,y\}$ of
vertices. Choose such a vertex $x$ that is contained in exactly
$r$ classes. Now $\sum_{i\colon x\in V_i} |V_i|\geq n+r-1$ since we count
every vertex at least once and $x$ exactly $r$
times. On the other hand, $\sum_{i\colon x\notin V_i} |V_i|\geq (l-r)k$
because the size of every class is at least $k$. Combining these
inequalities we obtain $\sigma=\sum_{i\colon x\in V_i} |V_i|+\sum_{i\colon x\notin V_i}
|V_i|\geq n+r-1+(l-r)k$, which proves the second part of the
lemma. $\square$
\end{biz}

Using the lemma, it is enough to show that either $rn\geq
n+(l-1)(k-1)$ or $n+r-1+(l-r)k\geq n+(l-1)(k-1)$ holds, since these
imply $\sigma\geq n+(l-1)(k-1)$ and therefore $m\geq n-(k-1)$ will
hold.

It is easy to see that
obtain $$rn\geq n+(l-1)(k-1)\, \Leftrightarrow\, l-1 \leq
\frac{(r-1)n}{k-1},$$ and
$$n+r-1+(l-r)k\geq n+(l-1)(k-1)\, \Leftrightarrow\, l-1 \geq(r-1)(k-1).$$
Therefore our claim, $m\geq n-(k-1)$, does not hold only if
$\frac{(r-1)n}{k-1}< l-1 < (r-1)(k-1)$. This implies the condition $n<(k-1)^2$ (we have
seen above that $r\neq 1$), so $n\geq (k-1)^2$ implies $m\geq
n-(k-1)$, which completes the proof. $\square$
\end{biz}

The lower bound of Theorem \ref{Tlower} is sharp since chains and
stars on $n$ vertices have exactly $n-(k-1)$ edges.

It is an interesting question if it is possible to eliminate the
condition $n\ge (k-1)^2$. Or otherwise, what are the sharp lower and
upper bounds for sizes of $k$-uniform counterexamples?

 We have partial
answers for these questions. Let us call a hypertree a counterexample
if it does not fulfil $m\geq n-(k-1)$.

\begin{allitas} \mbox{}
\begin{itemize}
\item For $k=3,4,5$, there are no $k$-uniform counterexamples, so
Theorem \ref{Tlower} holds without any condition for $n$.
\item For every $k\geq 6$, there is a $k$-uniform counterexample.
\end{itemize}
\end{allitas}

The proof of the first part is a fairly straightforward case analysis
using combinatorial arguments and relations between the parameters
$n,m,k,r$ and $l$ (see $r$ and $l$ in the proof of Theorem
\ref{Tlower}), so it is left to the reader.  The second part is
obtained from the following simple construction.

 Let
$V'=\{v_1, v_2, \ldots, v_{k-6}\}$, $V_x=\{x_1,x_2,x_3\}$,
$V_y=\{y_1,y_2,y_3\}$, $V_z=\{z_1,z_2,z_3\}$, $V=V_x\cup V_y\cup
V_z\cup V'$ and $\mathcal{E}=\{V_x\cup V'\cup V_y, V_x\cup V'\cup V_z,
V_y\cup V'\cup V_z\}$. Then the hypergraph $\mathcal{F}=(V,
\mathcal{E})$ is an edge-minimal hypertree that has $k+3$ vertices
and $3$ edges. Hence $m=3<4=n-(k-1)$, so $\mathcal{F}$ is a
counterexample if $V'$ exists, i.e. $k\geq 6$.

Let us see, how many vertices a counterexample may have for a fixed
$k$.
\begin{allitas} \mbox{}
\begin{itemize}
\item If $k\geq 6$, then there exist a $k$-uniform counterexample of
order $k+3$;
\item if $k\geq 6$ even, then there exist a $k$-uniform counterexample
of order $\frac{k(k-2)}{2}$;
\item if $k\geq 6$ odd, then there exist a $k$-uniform counterexample
of order at least $\frac{(k-1)(k-4)}{2}+1$;
\end{itemize}
\end{allitas}

The previous construction proves the first part. It can be easily seen
that it is a sharp lower bound since every hypertree with $2$ edges
must have $k+1$ vertices, thus a counterexample with $k+2$ vertices
cannot exist.

For the second statement, we use $c$ clusters of vertices, each of
size $\frac{k}{2}$, where $c$ depends on $k$, and each pair of
clusters form an edge.

Now, we have $n=\frac{ck}{2}$ vertices and $m=\binom{c}{2}$ edges. The
condition being a counterexample can be formulated
as $$\binom{c}{2}<\frac{ck}{2}-(k-1).$$ To maximize $n$ for a fixed
$k$, we must maximize the integer $c$ subject to the previous
constraint. $c$ attains its maximal value at $k-2$, and the hypergraph
(which is obviously chain-connected) obtained in this way has
$\frac{k(k-2)}{2}$ edges.

The third claim can be proven with a minor modification of the
previous proof. Take the construction above with $c$ clusters, each of
size $\frac{k-1}{2}$. Now add an extra vertex
to each edge, and maximize $c$ subject to the related, modified
constraint.

\section{Upper bounds for the edge number of hypertrees}\label{3}

In this section we prove an upper bound for the
edge number. It is also shown that the given bound is asymptotically
sharp in $3$-uniform case.

It is obvious that a $k$-uniform hypergraph has at most $\binom{n}{k}$
edges. In the case of hypertrees, the order of magnitude is one less.

\begin{tetel}\label{Tupper} If $\mathcal{F}=(V,\mathcal{E})$ is a
semicycle-free, $k$-uniform hypergraph with $n$ vertices and $m$
edges, then $m\leq \binom{n}{k-1}$.
\end{tetel}

\begin{biz} We give an injective function $\varphi \colon  \mathcal{E}
\rightarrow \binom{V}{k-1}$ below.

 Let us construct a
maximal chain in $\mathcal{F}$ (i.e. it cannot be continued with
more edges). Such a chain exists because $\mathcal{F}$ is
semicycle-free. Starting with any edge, we construct the chain edge by
edge in greedy way by appending a new edge to the end of the
chain. The semicycle-free property implies that every chain is
non-self-intersecting, hence the chain expands with a new vertex
(which is not used so far) in every step. Sooner or later the vertices
of $\mathcal{F}$ run out, and we cannot expand our chain anymore.

 Let us take the last edge of a maximal chain, and assign the
set consisting of the last $k-1$ vertices of the chain to this
edge. Obviously, this $(k-1)$-set is contained only in the edge to
which we assigned it, otherwise the chain could be continued.

 Now, let us delete this edge from $\mathcal{F}$, and repeat
the same process to the remaining hypergraph. The subhypergraphs
inherit the semicycle-free property, so we can execute the assignment
again. We continue deleting edges until the edges of $\mathcal{F}$ run
out, i.e. $\varphi$ is defined everywhere.

It remains to prove that $\varphi$ is injective. Assume indirectly
that we assign the same vertex set $T$ to two different edges $e_1$
and $e_2$. Without loss of generality, we may assume that during the
assigning process we have deleted $e_1$ earlier than $e_2$. So, $e_2$
still exists at the moment of deleting $e_1$, and both edges contains
the set $T$ because of the definition of $\varphi$. This is a
contradiction, thus $\varphi$ must be injective.

 Due to the
injectivity we have $|\varphi(\mathcal{E})|=m$, combining this with
$|\varphi(\mathcal{E})|\leq \binom{n}{k-1}$ we obtain the
theorem. $\square$
\end{biz}

\noindent \textbf{Remark:} Theorem \ref{Tupper} remains true even if
we only require $\mathcal{F}$ to be cycle-free because we can also
find a maximal chain in that case (otherwise there would be a
$k-1$-tuple of vertices, which appears at least twice as $k-1$
consecutive vertices of the chain).  The same proof also works in this
case.

In particular, the statement is true
also for hypertrees (since they are semicycle-free).

The following construction proves the asymptotic sharpness of the
upper bound of Theorem \ref{Tupper} in 3-uniform case.

Let $V=\{v_i\}_{i=1}^n$ be a set of vertices and
$\mathcal{F}=(V,\mathcal{E})$ be an arbitrary 3-uniform hypertree. Now let
$B(\mathcal{F})=(V\cup V',\mathcal{E}\cup\mathcal{E}')$ denote the
3-uniform hypergraph where $V'=\{0,1\}^n$,
\[\mathcal{E}'=\left\{\{v_i,u,w\}\colon  \text{
  \begin{minipage}[c]{0.5\linewidth}
    $v_i\in V$, $u,w\in V'$ and the $i$th bit is the first bit where
    $u$ and $w$ differ
  \end{minipage}
}\right\}.\]

\begin{tetel}\label{Thuen1} If $\mathcal{F}$ is a hypertree, then
$B(\mathcal{F})$ is also a hypertree.
\end{tetel}

\begin{biz}

\noindent $(1)$ Chain-connectedness:

 Any two vertices from
$V$ are connected by a chain because $\mathcal{F}$ is a hypertree,
and  all of its edges belong to $B(\mathcal{F})$, too.

For $u,w\in V'$ let $i$ denote the position of the first bit where
they differ. Then by definition $\{v_i,u,w\}\in\mathcal{E}'$, thus $u$
and $w$ are connected by a chain of length one in
$B(\mathcal{F})$.

 In case of $u\in V'$, $v_i\in V$ consider
a vertex $w\in V'$ whose first $i-1$ bits are the same as the first
$i-1$ bits of $u$, but its $i\text{th}$ bit differs from the
$i\text{th}$ bit of $u$. Such a $w$ certainly exists since the
preceding condition can be satisfied due to $1\leq i\leq n$. Then, by
definition, $\{v_i,u,w\}\in\mathcal{E}'$, so $u$ and $v_i$ are connected
with a chain of length one.


\noindent $(2)$ Semicycle-freeness:

 Assume indirectly that
$B(\mathcal{F})$ contains a semicycle $C$. If the first edge of $C$
belongs to $\mathcal{E}$, then $\mathcal{F}$ entirely contains $C$,
otherwise the first edge of $C$ that intersects $V'$, would include
exactly two vertices from $V$, but such an edge does not
exist. However, $\mathcal{F}$ cannot contain $C$ since $\mathcal{F}$
is semicycle-free, which leads to contradiction.

 Hence, we have the first edge of $C$ in $\mathcal{E}'$. Then
the second edge contains at least one vertex from $V'$, so this
edge also belongs to $\mathcal{E}'$. The same argument shows that
every edge of $C$ is in $\mathcal{E}'$. Let us notice that only one
edge can contain a pair $\{u,w\}\subseteq V'$ since the first place
where the bits of $u$ and $v$ differ, uniquely determines the third
vertex of the edge containing $\{u,w\}$, which is in $V$. So, if two
adjacent vertices of the semicycle are in $V'$, then these vertices
are the last two vertices of the semicycle since only one edge of $C$ can
contain them.

 Let $\{v_{i_1},u_1,u_2\}$ be the first edge of
$C$. If we write down the vertices of the semicycle in a sequence,
denoting the vertices from $V$ by $v_{i_k}$ and those from $V'$ by
$u_{k}$, there are 3 possible sequences:
\begin{itemize}
\item $v_{i_1} u_1 u_2\cdot$ : this cannot be a semicycle because it
would have less than 5 vertices,
\item $u_1 v_{i_1} u_2 u_3\cdot$ : the same argument works in this
case,
\item $u_1 u_2 v_{i_1} u_3 u_4\cdot$ : this cannot be a semicycle
since $u_1$ and $u_2$, $u_2$ and $u_3$, $u_3$ and $u_4$ differs in the
$i\text{th}$ bit, so $u_1$ and $u_4$ differs in the $i\text{th}$ bit,
but they have to be identical in the semicycle.
\end{itemize} Here, we have taken into consideration that every edge of
$C$ comes from $\mathcal{E}'$, and such an edge contains two
vertices from $V'$. The point at the end of the sequences indicates
the fact that the given sequence cannot be continued.

 We get
a contradiction again, so there is no semicycle in
$B(\mathcal{F})$,  implying that $B(\mathcal{F})$ is a
hypertree. $\square$
\end{biz}

Let us count the number of vertices and edges of the hypertree obtained
in this way. $|V\cup V'|=n+2^n$, $|\mathcal{E}\cup\mathcal{E}'|\geq
|\mathcal{E}'|=\binom{2^n}{2}$ since exactly one edge of
$\mathcal{E'}$ belongs to every pair of vertices of
$V'$. 
So, $$\frac{|\mathcal{E}(B(\mathcal{F}))|}{\binom{|V(B(\mathcal{F}))|}{2}}\geq
\frac{\binom{2^n}{2}}{\binom{n+2^n}{2}}\rightarrow 1\ \text{if
}n\rightarrow\infty.$$ This matches the bound in Theorem
\ref{Tupper}.

 Thus, if $\mathcal{F}_1, \mathcal{F}_2,
\ldots$ is a hypertree sequence for which $\lim_{n\rightarrow\infty}
|V(\mathcal{F}_n)|=\infty$ is true, then
$|\mathcal{E}(B(\mathcal{F}_n))|\sim\binom{|V(B(\mathcal{F}_n))|}{2}$.

If we want to get a hypertree with relatively few vertices and
high number of edges, then we can apply the operator $B()$ more times one
after the other.

This result verifies that the bound of Theorem \ref{Tupper} is
asymptotically sharp in 3-uniform case.  In $k$-uniform case, it is an
open question to find a similar construction.

\noindent \textbf{Remark:} Theorem \ref{Thuen1} implies that from
the viewpoint of the edge number, requiring cycle-free or
semicycle-free property does not mean a big difference in $3$-uniform case since the
cycle-free property guarantees that the number of edges is at most
$\binom{n}{2}$, and nor the stronger semicycle-freeness
can lower this bound.

\section{Upper bound for $l$-hypertrees}\label{4} 

We have seen the
upper bound for the edge number of a hypertree.  One can get better
bound under the assumption $\mathcal{F}$ contains no long chains.

\begin{tetel}\label{Tupper3} If $1\leq l\leq k$ and
$\mathcal{F}=(V,\mathcal{E})$ is a $k$-uniform
$l$-hypertree with $n$ vertices and $m$ edges, then $m\leq
\frac{1}{k-l+1}\binom{n}{k-1}$.
\end{tetel}

\begin{biz} The proof is similar to the proof of Theorem
\ref{Tupper}. We give $k-l+1$ injective functions,
$\varphi_1,\varphi_2,\ldots,\varphi_{k-l+1} \colon  \mathcal{E} \rightarrow
\binom{V}{k-1}$ for which
$\text{Im}\varphi_i\cap\text{Im}\varphi_{j}=\emptyset$, for all $1\leq
i<j\leq k-l+1$.

 Let $L$ be a maximal  chain in $\mathcal{F}$. Such a chain
exists because $\mathcal{F}$ is semicycle-free. Let $e$ be the last
edge of $L$. The length of this chain is at most $l\leq k$, therefore there
are at least $k-l+1$ vertices contained in the intersection of all edges of $L$. Let us
denote the set of these vertices by $U$. If $u\in U$ and some edge
$f\in \mathcal{E}$, $f\neq e$ covers $e\backslash \{u\}$, then $L$
could be continued by $f$ (because we have a freedom in ordering the
elements of $U$), which is a contradiction.

 Let $u_1,
\ldots, u_{k-l+1}\in U$ be distinct vertices and $\varphi_i
(e)=e\backslash\{u_i\}$, for all $i=1,\ldots, k-l+1$. Then delete $e$
and repeat the whole process for the remaining hypergraph. In the end,
$\varphi_i$s will be defined everywhere. The previous observation
shows that the conditions required for $\varphi_i$s are satisfied.

 So, $(k-l+1)|\mathcal{E}|=|\text{Im}
\varphi_1|+\ldots+|\text{Im}\varphi_{k-1}|=|\text{Im}
\varphi_1\cup\ldots\cup\text{Im}\varphi_{k-1}|\leq \binom{n}{k-1}$,
hence $|\mathcal{E}|\leq \frac{1}{k-l+1} \binom{n}{k-1}$. $\square$
\end{biz}

\noindent\textbf{Remark:} Semicycle-freeness is
a crucial premise of this proof, and there are no trivial extension to
the case when we exchange semicycle freeness by cycle freeness because
a maximal chain of length $l$ could be extended with an edge to form a
semicycle of length $l+1$.

 The most
important case is $l=2$, when the number of edges is at most
$\frac{1}{k-1}\binom{n}{k-1}$.


\section{Open problems}

There are many interesting open questions related to hypertrees such
that: ``What is the maximal number of edges in a $k$-uniform
edge-minimal hypertree of order $n$?'' or ``What is the minimal number
of edges in a $k$-uniform edge-maximal hypertree of order $n$?''.

Based on our research, we propose the following two
conjectures:

\begin{sejtes}\label{Smaxél} For every $k$-uniform edge-minimal
hypertree $\mathcal{F}=(V,\mathcal{E})$ on $n$ vertices
$|\mathcal{E}|\leq \frac{1}{k-1} \binom{n}{2}$ holds.
\end{sejtes}

\begin{sejtes}\label{Sminél} Every 3-uniform edge-maximal hypertree on
$n$ vertices has at least $\frac{1}{2}\binom{n}{2}-O(n)$ edges.
\end{sejtes}

It remained an open question whether or not the upper bound of the edge number of
$l$-hypertrees stated in Theorem \ref{Tupper3} is asymptotically sharp.


However, we can also modify the definition of edge-minimal hypertrees
a bit.  Instead of edge-minimal hypertrees, it is interesting to study
edge-minimal chain-connected hypergraphs or local hypertrees, which
can contain long semicycles, but no short ones. Similarly, we can study edge-maximal
semicycle-free hypergraphs instead of edge-maximal hypertrees.

  What
if we allow a chain to use an edge more times? Then our theorems and
definitions would alter more or less. In this case, one can show
forbidden substructures in edge-minimal hypertrees such as the complete hypergraph $K_{k+2}^{(k)}$.

A special kind of 3-uniform hypertrees is simultaneously edge-minimal and
edge-maximal 3-uniform hypertrees. This is a small subclass, and its
elements have nearly $\frac{1}{2}\binom{n}{2}$ edges by our
conjectures.

Another open problem is to give tight lower and upper bounds for the number
of edges in recursive hypertrees.

\section*{Acknowledgement}
The work reported in the paper has been developed in the framework of the
project "Talent care and cultivation in the scientific workshops of BME"
project. This project is supported by the grant T\'AMOP -
4.2.2.B-10/1--2010-0009.

The first author is partially supported by the Hungarian National
Research Fund (grant number NK 78439).


\begin{thebibliography}{00}
\bibitem{berge} C. Berge: Hypergraphs, North-Holland, Amsterdam, (1989)
 389--396.
\bibitem{looseham} D. K\"uhn, D. Osthus: Loose Hamilton Cycles in
3-uniform Hypergraphs of High Minimum Degree, \textit{Journal of
Combinatorial Theory, Series B} \textbf{96} (2006) 767--821.
\bibitem{lcycles} D. K\"uhn, R. Mycroft, D. Osthus: Hamilton $l$-cycles in
$k$-graphs, \textit{J. Combin. Theory Ser. A} \textbf{117} (2010) 
910--927.
\bibitem{hamchains} G. Y. Katona, H. Kierstead: Hamiltonian Chains in
Hypergraphs, \textit{Journal of Graph Theory} \textbf{30} (1999)
205--212.
\bibitem{loosecyc} H. H\`an, M. Schacht: Dirac-type Results for Loose
Hamilton Cycles in Uniform Hypergraphs, \textit{Journal of
Combinatorial Theory, Series B} \textbf{100} (2010) 332--346.
\bibitem{design} J. Demetrovics, G. OH. Katona, A. Sali: Design type
problems motivated by database theory, \textit{Journal of Statistical
Planning and Inference} \textbf{72} (1998) 149--164.
\bibitem{approx} V. R\"odl, A. Ruci\'nski, E. Szemer\'edi: An Approximate
Dirac-Type Theorem for $k$-Uniform Hypergraphs, \textit{Combinatorica}
\textbf{28} (2) (2008) 229--260.
\bibitem{szemeredi} V. R\"odl, A. Ruci\'nski, E. Szemer\'edi: Dirac-type
conditions for hamiltonian paths and cycles in 3-uniform hypergraphs,
\textit{Advances in Mathematics} \textbf{227} (2011) 1225--1299.
\end{thebibliography}
\end{document}